\documentclass[10pt,twoside]{article}
\usepackage{graphicx}
\usepackage{amsmath}
\usepackage{Latex-document}

\title{\bf Shock Waves}
\author{Tai-Ping Liu\thanks{Institute of Mathematics, Academia Sinica,
Taipei, China. Department of Mathematics, Stanford University, Stanford, CA 94305, USA. E-mail:
tpliu@math.sinica.edu.tw} \vspace*{-0.5cm}}
\date{\vspace{-8mm}}

\begin{document}

\maketitle

\thispagestyle{first} \setcounter{page}{185}

\begin{abstract}

\vskip 3mm

Shock wave theory was first studied for gas dynamics, for which
shocks appear as compression waves. A shock wave is characterized
as a sharp transition, even discontinuity in the flow. In fact,
shocks appear in many different physical situation and represent
strong nonlinearity of the physical processes.  Important
progresses have been made on shock wave theory in recent years.
We will survey the topics for which much more remain to be made.
These include the effects of reactions, dissipations and
relaxation, shock waves for interacting particles and Boltzmann
equation, and multi-dimensional gas flows.

\vskip 4.5mm

\end{abstract}

\vskip 12mm

\section{Introduction}

\vskip-5mm \hspace{5mm}

The most basic equations for shock wave theory are the systems of hyperbolic conservation laws
$$
u_t+\nabla_x \cdot f(u)=0, $$ where  $x\in R^m$ is the space
variables and $u\in R^n$ is the basic dependence variables. Such
a system represents basic physical model for which $u=u(x,t)$ is
the density of conserved physical quantities and the flux $f(u)$
is assumed to be a function of $u$. More complete system of
partial differential equations takes the form
$$u_t+\nabla_x \cdot
f(u)=\nabla_x\cdot (B(u,\varepsilon)\nabla_x u)+g(u,x,t), $$ with
$B(u,\varepsilon)$ the viscosity matrix and $\varepsilon$ the
viscosity parameters, and $g(u,x,t)$ the sources. Other
evolutionary equations which carry shock waves include the
interacting particles system, Boltzmann equation, and discrete
systems. Discrete systems can appear as difference approximations
to hyperbolic conservation laws. In all these systems, shock
waves yield rich phenomena and also present serious mathematical
difficulties due to their strong nonlinear character.

\section{Hyperbolic conservation laws}

\vskip-5mm \hspace{5mm}

Much has been done for hyperbolic conservation laws in one space
dimension
$$u_t+f(u)_x=0,\ x\in R^1,
$$
see \cite{smo}, \cite{daf}, \cite{lium}, and the article by Bressan in this
volume. Because the solutions in general contain discontinuous
shock waves, the system provides impetus for the introduction of
new ideas, such as the Glimm functional,  and is a good testing
ground for new techniques, such as the theory of compensated
compactness, in nonlinear analysis.

\vspace*{-2mm}

\section{Viscous conservation laws}

\vskip-5mm \hspace{5mm}

Physical models of the form of viscous conservation laws are not uniformly parabolic, but hyperbolic-parabolic.
Basic study of the dissipation of solutions for such a system has been done using the energy method, see
\cite{kaw}. Study of nonlinear waves for these systems has been initiated,  \cite{liu},  \cite{liuzeng}, however,
much more remains to be done. The difficulty lies in the nonlinear couplings due to both the nonlinearity of the
flux $f(u)$, which is the topic of consideration for hyperbolic conservation laws, as well as that of the
viscosity matrix $B(u,\varepsilon)$. For instance, the study of zero dissipation limit $\varepsilon \rightarrow
0$, see Bressan's article, has been done only for the artificial viscosity matrix $B(u,\varepsilon)=\varepsilon
I$.

\vspace*{-2mm}

\section{Conservation laws with sources}

\vskip-5mm \hspace{5mm}

Sources added to conservation laws may represents geometric effects, chemical reactions, or relaxation effects.
Thus there should be no unified theory for it. When the source represents the geometric effects, such the
multi-dimensional spherical waves , hyperbolic conservation laws takes the form
$$u_t+f(u)_r=\frac{m-1}{r}h(u),\ r^2=\sum_{i=1}^m (x_i)^2.$$
There has stablizing and destablizing effects, such as in the
nozzle flows,  \cite{lien}. The chemical effects occur in the
combustions. There is complicated, still mostly not understood,
phenomena on the rich behaviour of combustions. One interesting
problem is the transition from the detonations to deflagrations,
where the combined effects of dissipation, compression and
chemical energy gives rise to new wave behaviour. Viscous effects
are important on the qualitative behaviour of nonlinear waves
when the hyperbolic system is not strictly hyperbolic, see
 \cite{liume}. Relaxation, such as for the kinetic models and thermal
non-equilibrium in general, is interesting because of the rich
coupling of dissipation, dispersion and hyperbolicity,  \cite{cll}.

\vspace*{-2mm}

\section{Discrete conservation laws}

\vskip-5mm \hspace{5mm}

Conservative finite differences to the hyperbolic conservation
laws, in one space dimension, take the form:
$$u^{n+1}(x)-u^n(x)=\frac{\Delta t}{\Delta x}(F[u^n](x+\frac{\Delta x}{2})-F[u^n](x-\frac{\Delta x}{2})).$$
It has been shown for a class of two conservation laws and dissipative schemes,
such as Lax-Friedrichs and Godunov scheme that the numerical solutions converge to the exact solutions of the
conservation laws,  \cite{ding}. On the other hand, qualitative studies on the nonlinear waves for difference
schemes indicate rich behaviour. In particular, the shock waves depend sensitively on the its C-F-L speeds,
\cite{se}, \cite{ly1}.

\vspace*{-2mm}

\section{Multi-dimensional gas flows}

\vskip-5mm \hspace{5mm}

Shock wave theory originated from the study of the Euler
equations in gas dynamics. The classical book  \cite{cf} is still
important and mostly updated on multi-dimensional gas flows.
Because of its great difficulty, the study of multi-dimensional
shocks concentrate on flows with certain self-similarity
property. One such problem is the Riemann problem, with initial
value consisting of finite many constant states. In that case,
the solutions are function of $x/t$, not general function of
$(x,t)$, see  \cite{zhang}. See also  \cite{zheng} for other self-similar
solutions. However, unlike single space case, multi-dimensional
Riemann solutions do not represent general scattering data, and
are quite difficult to study. It is more feasible to consider
flows with shocks and solid boundary, e.g.  \cite{chen}  \cite{lien}.

\vspace*{-2mm}

\section{Boltzmann equation}

\vskip-5mm \hspace{5mm}

The Boltzmann equation $$f_t+\xi\nabla_x f=Q(f,f)$$ contains much
more information than the gas dynamics equations. Nevertheless,
the shock waves for all these equations have the same
Rankine-Hugoniot relation at the far states. The difference is on
the transition layer. There is beginning an effort to make use of
the techniques for the conservation laws to study the Boltzmann
shocks,  \cite{ly}. This is a line of research quite different from the
intensive current efforts on the incompressible limits of the
Boltzmann equation.

\vspace*{-2mm}

\section{Interacting particle systems}

\vskip-5mm \hspace{5mm}

Interacting particle systems is even closer to the first physical
principles than the Boltzmann equation. There is the
long-standing problem, the Zermelo paradox, in passing from the
reversible particle systems to the irreversible systems such as
the Boltzmann equation, and the Euler equations with shocks. This
is fine for particle system with random noises. However, except
for scalar models, so far the derivation of Euler equations from
the particle systems has been done only for solutions with no
shocks,  \cite{v}.

\end{document}